\documentclass[12pt, oneside]{article} 

\usepackage[utf8]{inputenc}
\usepackage[T1]{fontenc}

\usepackage{amsmath}
\usepackage{amssymb} 
\usepackage{mathrsfs} 

\usepackage{amsthm} 
\theoremstyle{plain}
\newtheorem{theorem}{Theorem} 
\newtheorem{definition}[theorem]{Definition}

\newtheorem{lemma}[theorem]{Lemma}

\usepackage{algorithm}
\usepackage[noend]{algpseudocode}
\makeatletter
\def\BState{\State\hskip-\ALG@thistlm}
\makeatother

\usepackage{bm} 
\usepackage{bbm} 

\usepackage{graphicx} 
\usepackage{subcaption}

\usepackage{tikz} 
\usetikzlibrary{positioning ,shadows,matrix}

\usepackage[english]{babel} 

\usepackage{hyperref} 
\hypersetup{
    pdfpagemode={UseOutlines},
    bookmarksopen,
    pdfstartview={FitH},
    colorlinks,
    linkcolor={blue}, 
    citecolor={blue}, 
    urlcolor={blue} 
}

\usepackage{natbib} 
\bibpunct{\!(}{)}{;}{a}{,}{,}

\usepackage{color} 

\makeatletter  
\@ifundefined{@currsizeindex} 
  {\RequirePackage{relsize}\let\dolarger\relsize} 
  {\def\dolarger#1{\larger[#1]}} 

\newcommand*\@@bigtimes[2]{\vphantom{\prod} 
  \vcenter{\hbox{\dolarger{4}$\m@th#1\mkern-2mu\times\mkern-2mu$}}} 
\newcommand*\bigtimes{\mathop{\mathpalette\@@bigtimes\relax}\displaylimits} 
\makeatother

\def\A{\mathbb{A}}\def\N{\mathbb{N}}\def\P{\mathbb{P}}
\def\Z{\mathbb{Z}}
\def\1{\mathbbm{1}}

\def\eb{\bold e}\def\Mb{\bold M}\def\mb{\bold m}

\newcommand{\eq}[1]{ 
\begin{equation}
\begin{split}
#1
\end{split}
\end{equation}
} 

\newcommand{\eqstar}[1]{
\begin{equation*}
\begin{split}
#1
\end{split}
\end{equation*}
} 

\usepackage{xspace}
\def\ESF{Ewens sampling formula\xspace}
\def\ESFs{Ewens sampling formulae\xspace}
\def\PSF{Pitman sampling formula\xspace}
\def\PSFs{Pitman sampling formulae\xspace}

\usepackage{authblk}

\title{\bf A reversible allelic partition process and Pitman sampling formula}
\author{\normalsize Matteo Giordano$^{*}$, Pierpaolo De Blasi$^{\dagger,\star}$ and Matteo Ruggiero$^{\dagger,\star}$}
\affil{\normalsize\emph{University of Cambridge$^{*}$, University of Torino$^{\dagger}$ and Collegio Carlo Alberto$^{\star}$}}

\allowdisplaybreaks

\begin{document}
\maketitle

\abstract{
We introduce a continuous-time Markov chain describing dynamic allelic partitions which extends the branching process construction of the \PSF in \cite{Pitman2006} and the classical birth-and-death  process with immigration studied in \cite{Karlin1967}, in turn related to the celebrated \ESF. A biological basis for the scheme is provided in terms of a population of individuals grouped into families, that evolves according to a sequence of births, deaths and immigrations. We investigate the asymptotic behaviour of the chain and show that, as opposed to the birth-and-death process with immigration, this construction maintains in the temporal limit the mutual dependence among the multiplicities.
When the death rate exceeds the birth rate the system is shown to have a reversible distribution, identified as a mixture of \PSFs, with negative binomial mixing distribution on the population size. The population therefore converges to a stationary random configuration, characterised by a finite number of families and individuals. 
}

\normalsize
\section{Allelic partitions and sampling formulae}
\label{Sec:Intro}

A partition of $n\in \N$ is an unordered collection $\pi=\{n_1,\dots,n_k\}$ of $k\le n$ positive integers whose sum equals $n$. A common equivalent way to describe $\pi$ is by means of the so-called \emph{allelic partition}, which groups the partition sets by size. We denote by $\mb=(m_1,\dots,m_n)\in\Z^n_+$ the associated vector of multiplicities, where $m_{i}$ 
counts the number of repetitions of $i$ in $\pi$, and let 
$$
\A_n=\bigg\{\mb\in\Z^n_+,\ \sum_{i=1}^nim_i=n\bigg\}
$$
be the finite set of all allelic partitions of $n$.
%
%
It will be useful to embed $\A_n$ into $\Z^\infty_+$ by considering infinite vectors $\mb=(m_1,\dots,m_n,0,0,\dots)$, defining also
$$
s(\mb):=\sum_{i\ge1}im_i,\quad  
k(\mb):=\sum_{i\ge1} m_i,
$$
which are respectively the number of items and groups (or positive multiplicities) in $\mb$, and
$$
\A:=\bigcup_{n\ge1}\A_n=\{\mb\in\Z^\infty_+, \ s(\mb)<\infty\},
$$
the countable set of all allelic partitions.
%
%
A random allelic partition $\Mb=(M_i)_{i\ge1}$ can then be defined as a random variable taking values in $\A$. 
%
%


The literature on the distributional properties of random allelic partition is rich and well established. 
\cite{Ewens1972} derived the distribution of the random allelic partition induced by a sample of $n$ genes drawn from a selectively neutral population at equilibrium (see also \citealp{Karlin1972}). This is described by the celebrated \emph{\ESF} (ESF), which assigns probability 
\eq{\label{ESF}
\textnormal{ESF}^\theta_n(\mb)=\frac{n!}{\theta_{(n)}}\prod_{i\ge1}\left(\frac{\theta}{i}\right)^{m_i}\frac{1}{m_i!}\1_{\{s(\mb)=n\}},
}
to the configuration $\mb\in \A$ with $m_i$ alleles appearing exactly $i$ times for each $i\ge1$. Here $\theta>0$ and $\theta_{(n)}=\theta(\theta+1)\dots(\theta+n-1)$ is the ascending factorial. The impact of the ESF has been significant in a number of fields beyond population genetics. For example, \cite{Antoniak1974} derived it independently from a  sample from a Dirichlet process, a cornerstone nonparametric prior distribution in Bayesian nonparametric statistics; \cite{Hoppe1984} recovered the formula as the marginal distribution of a Markov chain generated by a P\'olya-like urn model comprising an infinite number of colors, analogous to that introduced in \cite{Blackwell1973} and described equivalently in \cite{Aldous1985} with the famous metaphor of the \emph{Chinese restaurant process}. See \cite{johnson1997multivariate} and \cite{Crane2016} for reviews and \cite{Feng2010} for connections with population genetics. 

\cite{Pitman1995} introduced a two-parameter generalization of (\ref{ESF}), often referred to as the \emph{\PSF} (PSF), whereby for parameters $0\le\alpha<1$ and $ \theta>-\alpha$, the probability assigned to a random allelic partition $\mb\in \A$ is
\eq{\label{PSF}
\textnormal{PSF}_n^{\alpha,\theta}(\mb)=\frac{n!}{\theta_{(n)}}\prod_{i=1}^{k(\mb)}(\theta+i\alpha)\prod_{j\ge1}\left[\frac{(1-\alpha)^{j-1}}{j!}\right]^{m_j}\frac{1}{m_j!}\1_{\{s(\mb)=n\}}.
}
For $\alpha=0$, equation (\ref{PSF}) immediately recovers (\ref{ESF}), so that $\textnormal{PSF}_n^{0,\theta}\equiv\textnormal{ESF}_n^\theta$. For $0<\alpha<1$, the distribution arises for example in the study of stable processes with index $\alpha$, the case $\alpha=1/2$ being related to the zeros of Brownian motion (see also \citealp{Pitman1997a}), and in the study of the partition structures induced by a sample from a two-parameter Poisson--Dirichlet distribution \citep{Perman1992,Perman1993,Pitman1997b}.
The corresponding generalization of Hoppe's urn is defined in \cite{Pitman1995,pitman1996b} as a Markov chain of allelic partitions with initial state $\eb_0$ and transition probabilities
\eq{\label{PitTransProb}
p(\mb'|\mb)\propto
\begin{cases}
\theta+\alpha k(\mb), & \mb'=\mb+\eb_1,\\
(i-\alpha)m_i,& \mb'=\mb-\eb_i+\eb_{i+1}, \ i\ge1,\\
0, & \textnormal{else},
\end{cases}
}
where $\eb_i=(\delta_{ij})_{j\ge1}, \ i\ge0$, and the normalising constant is given by $\theta+s(\mb)$. The first transition in (\ref{PitTransProb}) can be understood as the addition to the gene sample of a new gene of a previously unobserved allele, while the transition from $\mb$ to  $\mb-\eb_i+\eb_{i+1}$ corresponds to adding a gene whose allele had been previously observed $i$ times. For $\alpha=0$, (\ref{PitTransProb}) coincides with Hoppe's urn, while for any $0\le\alpha<1$ the marginal distribution of the state of the chain after $n$ transitions is given by (\ref{PSF}). 

Among other relevant contributions to the related theory, the asymptotic behavior of the number  $K^{\alpha,\theta}_n$ of different alleles was derived by \cite{Korwar1973} in the case $\alpha=0$, for which
\eq{\label{KnConv1}
\lim_{n\to\infty}\frac{K^{0,\theta}_n}{\log n}=\theta\  \textnormal{a.s.},
} 
and by \cite{Pitman1997a} for $0<\alpha<1$, where we have that
\eq{\label{KnConv2}
\lim_{n\to\infty}\frac{K^{\alpha,\theta}_n}{\log n}=S^{\alpha,\theta} \ \textnormal{a.s.},
}
$S^{\alpha,\theta}$ being an absolutely continuous random variable on $(0,\infty)$ whose law is related to the Mittag--Leffler distribution. Concerning the PSF, we further refer the reader to \cite{Kerov2006,Pitman1996a} for various characterizations thereof, \cite{James2008} for connections with Bayesian nonparametric statistics and again \cite{Feng2010} for connections with population genetics. Note that  similar asymptotic results to the above are available for Gibbs-type models \citep{gnedin2006exchangeable,de2015gibbs}, which include the one- and two-parameter models recalled above, as well as a model in \cite{gnedin2010} which produces a finite but random number of families in the limit.

In this paper we are interested in the connection between Ewens--Pitman sampling formulae and the theory of stochastic processes for population growth. In particular, we construct a birth-and-death process with immigration whose dynamics modifies the one arising as after \eqref{PitTransProb} by allowing the removal of items from the system; and we show, in a particular regime, that the process is reversible with respect to a certain mixture of \PSFs. The rest of the paper is organised as follows. In Section \ref{sec: Main contribution} we present our contribution and its connections with past related work. Section \ref{Sec:DefMod} describes in detail the birth-and-death population model with immigration, and how it can be constructed via branching processes. Finally, Section \ref{Sec:Reversibility} analyses the reversible regime, identifying explicitly the reversible distribution and its connections with \PSFs.



\section{Main contribution and related work}\label{sec: Main contribution}

The seminal work of \cite{Karlin1967} describes a population wherein new families are initiated at a sequence of random times generated by a stochastic process $I=\{I(t),t\ge0\}$, which can be thought of as immigration events,
%
%
and then evolve, independently of one another, according to the law of a common continuous-time Markov chain on $\Z_+$, whose infinitesimal rates are denoted $q_{ij}$.
As pointed out by \cite{Tavare1989}, if new families are interpreted as novel mutant alleles of a given gene, then the scheme may be regarded as a version of the infinitely-many-neutral-alleles model and can provide, under specific choices of its probabilistic components, a generating mechanism for the sampling formulae (\ref{ESF}) and (\ref{PSF}). In particular, let $I$ be a pure-birth processes started at $0$ with birth rates 
$$
\lim_{h\to0}\frac{1}{h}\P(I(t+h)-I(t)=1|I(t)=k)=\theta+\alpha k, \quad  k\ge0,
$$
and assume that the process describing the evolution of each family is started at 0 and has rates
\eq{\label{FengBirthRate2}
q_{i,i+1}=i-\alpha,\ i\ge1.
}
If now $M_i(t)$ is the number of families of size $i$ at time $t$, and $\Mb(t)=(M_i(t))_{i\ge1}$ is the induced allelic partition, then $\Mb=\{\Mb(t),t\ge0\}$ defines a continuous-time Markov chain on $\A$ whose embedded jump chain has transition probabilities (\ref{PitTransProb}) and marginal distributions described by the PSF; see \cite{Feng1998} and \cite{Pitman2006}. When $\alpha=0$, this construction reduces to the birth process with immigration considered in \cite{Tavare1987}, which provides the corresponding embedding of Hoppe's urn leading to the ESF (cf.~Section 2.7.2 in \citealp{Feng2010}).

In the above schemes, each transition entails the addition of one individual to the population, and both the population size and the number of families diverge almost surely as time increases. This is a convenient mathematical simplification but may be undesirable or result in a lack of flexibility in different contexts, and is arguably unrealistic when modelling the evolution of populations. Therefore, it can be of interest to account for the death of individuals and to study the related implications on the family structure.

Motivated by the above considerations, we consider here  a modified version of the population model underlying the PSF, by assuming that the individuals are endowed with independent random exponential lifetimes of parameter $\mu>0$, with $\mu=0$ recovering the original construction, thereby introducing uniform death events in the scheme. See Definition \ref{Def:ThreeParProc} below for the details. By leaving unchanged the other rules governing the evolution of the family structure, we thus obtain more flexible dynamics whereby families are started, fluctuate in size and possibly become extinct with the passage of time, resulting in a population with varying size and number of families. Note that, because of exchangeability, the PSF is consistent with respect to uniform deletion, hence a death event induces a random partition still following the PSF. For consistency with respect to deletion of an entire family, a property known in the literature as regeneration, see \cite{GnedinPitman2005}.


For $\alpha=0$, our construction recovers the detailed version of the birth-and-death process with immigration (BDI) in \cite{Karlin1967}, object of thorough study in \cite{Kendall1975} and \cite{Tavare1989}. In this case, the new families are started at the times of a Poisson point process of intensity $\theta$, and each family subsequently evolves according to an independent linear birth-and-death process. Theorem 2.1 in \cite{Karlin1967} derives the marginal distribution of the induced allelic partition $\Mb(t)$, which is that of a sequence of independent Poisson random variables given by
\eq{\label{TwoParMargDistr}
\Mb(t)\sim\prod_{i\ge1}\textnormal{Po}(\theta b_t^i/i), \ 
b_t=
\begin{cases}
[e^{(1-\mu)t}-1]/[e^{(1-\mu)t}-\mu], & \mu\neq1\\
t/(1+t), & \mu=1.
\end{cases}
%
%
}
Anticipated by earlier results in \cite{Watterson1974} and \cite{Kendall1975},  \cite{Tavare1989} reformulated the above distribution as a mixture of \ESFs in (\ref{ESF}) with a negative binomial mixing distribution on $n$. Specifically, letting $S(t)$ being the number of alive individuals at time $t$, the population size process $S=\{S(t), t\ge0\}$ defines by construction a BDI process, with negative binomial marginal distribution
 \eq{\label{PopSizeMargDistr}
S(t)\sim \textnormal{N-Bin}(\theta,b_t),\ 
\P(S(t)=n)=\frac{\theta_{(n)}}{n!}(1-b_t)^\theta b_t^n, \ n=0,1,\dots
}
Then,  by combining (\ref{TwoParMargDistr}) and (\ref{PopSizeMargDistr}), the ESF is recovered as the conditional distribution of the actual allelic partition given the population size, i.e., 
\eq{\label{TwoParCondDistr}
\P(\Mb(t)=\mb|S(t)=n)=\textnormal{ESF}_n^\theta(\mb), \ \mb\in \A, \ n\ge1.
}

Within the study of the BDI process, a great attention has been dedicated to investigating the long-run behavior of the model. By taking the limit as $t\to\infty$ in (\ref{TwoParMargDistr}), we have the convergence in law
\eq{\label{TwoParConv}
\Mb(t)\overset{d}{\to}(X_1,X_2,\dots), \ t\to\infty,
}
where, for $0<\mu\le1$, $X_i\overset{\textnormal{ind}}{\sim}\textnormal{Po}(\theta/i)$, while for $\mu>1$, $X_i\overset{\textnormal{ind}}{\sim}\textnormal{Po}(\theta\mu^{-i}/i)$. In the first case  the death rate does not exceed the birth rate, and the asymptotic regime corresponds to the weak limit as $n\to\infty$ of the Ewens partition structure (\ref{ESF}) derived in \cite{Arratia1992}, characterized by an infinite sample size and a logarithmic growth of the  underlying number of groups.
%
%
Instead, for $\mu>1$, the death events occur at a faster rate than births, causing each family to eventually become extinct almost surely (see \citealp{Kendall1975}) and preventing the indefinite growth observed previously. In fact, in the second regime the population size process $S$ is easily seen in Lemma \ref{Lem:SReversibility} below to have $\textnormal{N-Bin}(\theta,\mu^{-1})$ reversible distribution. In turn, \cite{Kendall1975} showed the reversibility of the induced allelic partition process, and shed light on the representation of the reversible distribution as an analogous mixture as that arising from (\ref{TwoParCondDistr}), characterized by a $\textnormal{N-Bin}(\theta,\mu^{-1})$ mixing on the population size. 


In this paper, 
we encode the description of the evolving family structure directly in terms of the induced allelic partition process, and then focus on investigating the long-run behavior of the model. Because of the more involved rules governing the formation of new families,
%
%
the sharp distributional result derived in \cite{Karlin1967} is not accessible in our construction. Nonetheless, by observing that the overall dynamics of the total population size is left unchanged, we show that for $\mu>1$ the model is reversible, and identify the reversible distribution as a mixture of Pitman sampling formulae with respect to the same negative binomial mixing on the population size appearing for $\alpha=0$. For $0<\alpha<1$, the mixture can be written in closed form as
\eq{\label{ThreeParRevDistr}
\pi
(\mb)=C(\theta/\alpha)_{(k(\mb))}\prod_{i\ge1}\mathrm{Po}(m_i;\alpha_i\mu^{-i}),\ \alpha_i=\frac{\alpha(1-\alpha)_{(i-1)}}{i!}, \ \mb\in \A,
}
where 
\eq{\label{CNormConst}
C
=e^{1-(1-1/\mu)^{\alpha}}(1-1/\mu)^\theta.
}
%
%
Thus, in the reversible regime, the system will converge to a stationary random configuration that features almost surely finitely many families and individuals. However, as opposed to the previous model, the mutual dependence among the multiplicities is seen in (\ref{ThreeParRevDistr}) to be preserved in the limit.

\section{A birth-and-death process with immigration}
\label{Sec:DefMod}



\subsection{Definition}

\begin{definition}\label{Def:ThreeParProc}
For $0\le\alpha<1$, $\theta>-\alpha$ and $\mu>0$, 
let $\Mb=\{\Mb(t), \ t\ge0\}, \ \Mb(t)=(M_i(t))_{i\ge1}$, be a 
%
%
continuous-time Markov chain  with state space $\A$, initial state $\Mb(0)=\eb_0$ and temporally homogeneous transition rates
\eq{\label{ThreeParTransRate}
q(\mb'|\mb)&=
\begin{cases}
\theta+\alpha k(\mb),& \mb'=\mb+\eb_1,\\
(i-\alpha)m_i,&  \mb'=\mb-\eb_i+\eb_{i+1},\ i\ge1,\\
\mu im_i,& \mb'=\mb-\eb_i+\eb_{i-1},\ i\ge1,\\
-\theta-(1+\mu)s(\mb), & \mb'=\mb,\\
0,& \mathrm{else}.
\end{cases}\\
}
\end{definition}

%
%

We notice that the first two lines in (\ref{ThreeParTransRate}) coincide precisely with the transitions in the sequential construction by \cite{Pitman1995} displayed in (\ref{PitTransProb}). These describe the initiation of a new family of size one and the addition of one individual to a pre-existing family of size $i$ respectively. The third transition instead represents, for $i\ge2$, the death of a member of a family of size $i$, while for $i=1$ it corresponds to the extinction of a family comprising a single individual. The choice of the initial state $\Mb(0)=(0,0,\dots)$ is rather common across the literature, and can be interpreted as the population being initially empty, as in a territory about to be colonised.


The infinitesimal generator corresponding to the transition rates in (\ref{ThreeParTransRate}) 
are easily seen to be \textit{stable}, as $q(\mb|\mb)$ is finite for all $\mb\in\A$, and \textit{conservative}, since each row has zero sum (cf.~\citealp{Norris1997}). Furthermore, by adapting the argument in Section 4 in \cite{Kendall1975}, we can deduce from the properties of the associated population size process that $\Mb$ is also \textit{regular}, i.e.~it performs almost surely only a finite number of jumps in every finite time interval. In particular, for any $t\ge0$, let
\eq{\label{size process}
S(t)&:= s(\Mb(t))=\sum_{i\ge1}iM_i(t)
}
be the random number of items underlying the allelic partition $\Mb(t)$. Then we have from Definition \ref{Def:ThreeParProc} that $S=\{S(t),t \ge0\}$ defines a continuous-time Markov chain on $\Z_+$, with initial state $S(0)=0$ and transition rates
\eq{\label{STransRate}
r(n'|n)=
\begin{cases}
\theta+n, & n'=n+1,\\
\mu n, & n'=n-1,\\
-\theta-(1+\mu)n,&n'=n,\\
0,&\textrm{otherwise}.
\end{cases}
}
Indeed, by properly aggregating the rates in (\ref{ThreeParTransRate}), we see that conditionally given the population size $S(t)=n$ (that is to say $\Mb(t)=\mb$ for any $\mb\in\A$ with $s(\mb)=n$) one of the following events occur. A new individual is introduced in the population (either initiating a new family or joining a pre-existing one) increasing its size by one unit, and this happens with overall rate
$$
\theta+\alpha k(\mb)+\sum_{i\ge1}(i-\alpha)m_i
=\theta+n.
$$ 
Alternatively, the death of one of the $n$ individuals occurs, decreasing the population size by one unit, which happens at the total rate 
$$
\sum_{i\ge1}\mu im_i=\mu n.
$$
Thus, irrespective of the parameter $\alpha$, $S$ defines a classical BDI process, with unit birth rate, death rate $\mu$ and immigration rate $\theta$. The law of $S$ is seen in (\ref{STransRate}) to be insensitive  to the choice of $\alpha$, so that we conclude that the more refined rules governing the formation of new families do not impact the dynamics of the total population size, but rather influence the detailed probabilistic fluctuations in the family configuration. 

Finally, we notice that $\Mb$ and $S$ share the same random jump time, as each transition of the allelic partition process involves either the introduction or the death of an individual. Therefore, as the generator of the birth-and-death process with immigration is known to be regular (see, e.g., Section 2 in \cite{Kendall1975}), we can directly conclude about the regularity of $\Mb$, which accordingly identifies (jointly with the initial state) a well-defined chain on $\A$.

%
%


\subsection{Branching process construction}

A biological basis for the three-parameter process in Definition \ref{Def:ThreeParProc} can be provided in terms of a population of individuals grouped into families evolving according to a sequence of births, deaths and immigrations, by introducing uniform death events in the models studied by \cite{Feng1998} and \cite{Pitman2006}.

First assume that $\theta>0$, and consider a population whose members have independent exponential lifetimes with parameter $\mu$, and autonomously reproduce according to independent unit-rate Poisson processes (taken without loss of generality, upon appropriate time rescaling). At a birth, any offspring produced by the oldest member of each family will either join its parent's family, with probability $1-\alpha$, or start a new family, with probability $\alpha$. All other family members give birth to individuals that will remain within their parent's family. Finally, an independent Poisson process of rate $\theta$ is assumed to bring immigrants into the population that, at the time of entry, will initiate a new family each.

For $-\alpha<\theta\le0$, the above rules can be modified by assuming that no immigrant can enter the population, but rather that the overall oldest member in the population produces at  rate $1-\alpha$ offspring that will join its family, and at rate $\alpha+\theta$ individuals that will initiate new families. 

Thus, by exploiting basic properties of independent Poisson processes (e.g., a \textit{Poissonization} argument as in \citealp{Athreya1968}), it follows that, if we encode the family structure into the corresponding sequence of multiplicities, the dynamics of the above population  generates a continuous-time Markov chain of allelic partitions, whose transition rates recover precisely those of the process $\Mb$ in (\ref{ThreeParTransRate}). In particular, each family evolves according to an independent birth-and-death process with initial state 1 and transition rates
$$
q_{i,i+1}=i-\alpha,\quad  q_{i,i-1}=\mu i, \quad  i\ge1
$$
directly generalizing (\ref{FengBirthRate2}). 
On the other hand, the introduction of deaths changes substantially the probabilistic laws that govern the input process $I=\{I(t), \ t\ge0\}$ (counting the number of new families formed up to time $t$). In particular, as opposed to \cite{Karlin1967}, here $I$ does not define a Markov process, since the probability of introducing a new family at each transition depends on the actual number of families $k(\mathbf M(t))$ which, due to families possibly becoming extinct, is in general different from $I(t)$. This in turn represents a serious obstruction in deriving distributional results for quantities of interest related to the model, in particular for the total number of families. Still, as resulting from (\ref{STransRate}), the dynamics of the total population size is tractable, and this fact will be exploited in the next section to study a reversibility regime for the process.

\section{The reversible regime}
\label{Sec:Reversibility}


As seen above, the population size process $S$ defines a classical BDI process. Hence, from the transition rates in (\ref{STransRate}) we immediately deduce that $S$ is irreducible, as each pair of states are mutually accessible. Moreover, the process is known to be \textit{transient} for $0<\mu<1$, \textit{null recurrent} for $\mu=1$ and \textit{positive recurrent} for $\mu>1$ (cf.~Section 2 in \citealp{Kendall1975}). In the latter case, $S$ is  \textit{ergodic}, 
%
%
implying the convergence in law
$$
\lim_{t\to\infty}\P(S(t)=n|S(0)=0)=\lambda(n), \ \forall n\in\Z_+
$$ 
for a unique \textit{stationary distribution} $\lambda=(\lambda(n))_{n\ge0}$ on $\Z_+$ satisfying
$$
\sum_{n\in\Z_+}\lambda(n)r(n'|n)=0, \ \forall n'\in\Z_+.
$$
The above sequence of equations can be solved recursively to obtain
\eq{\label{SRevDistr}
\lambda(n)=\frac{\theta_{(n)}}{n!}\mu^{-n}(1-1/\mu)^\theta, \ n\ge0,
}
%
%
so that $\lambda$ results to be a negative binomial distribution $\textrm{N-bin}(\theta,\mu^{-1})$ with parameters $\theta$ and $\mu^{-1}$. The following lemma shows that, under the condition $\mu>1$, the population size process in fact enjoys the stronger property of being reversible with respect to $\lambda$.

\begin{lemma}\label{Lem:SReversibility}

Let $S$ be the population size process associated to $\Mb$ in Definition \ref{Def:ThreeParProc} through \eqref{size process}. Then, if $\mu>1$, $S$ is reversible with respect to its stationary distribution $\lambda$ in \eqref{SRevDistr}.

\proof

As $S$ is regular and irreducible, the reversibility of the process is equivalent to the {detailed balance condition}
\eq{\label{SBalCond}
\lambda(n)r(n'|n)=\lambda(n')r(n|n'),\ n\neq n'.
}
In view of (\ref{STransRate}), we only need to show (\ref{SBalCond}) for the cases $n\ge0, \ n'=n+1$, and $n\ge1, \ n'=n-1$, as the transition rates corresponding to all other transitions are equal to 0.
For $n\ge0$, $n'=n+1$, the left hand side of (\ref{SBalCond}) becomes
\eqstar{
\lambda(n)r(n+1|n)&=\dfrac{\theta_{(n)}}{n!}\mu^{-n}(1-1/\mu)^\theta (\theta+n)\\
&=\dfrac{\theta_{(n+1)}}{(n+1)!}\mu^{-(n+1)}(1-1/\mu)^\theta \mu(n+1)\\
&=\lambda(n+1)r(n|n+1);
}
and similarly, for $n\ge1$, $n'=n-1$ we have
\eqstar{
\lambda(n)r(n-1|n)&=\dfrac{\theta_{(n)}}{n!}\mu^{-n}(1-1/\mu)^\theta\mu n \\
&=\dfrac{\theta_{(n-1)}}{(n-1)!}\mu^{-(n-1)}(1-1/\mu)^\theta (\theta+n) \\
&=\lambda(n-1)r(n|n-1).
}

\endproof

\end{lemma}


We now move on to  study the properties of $\Mb$. We first notice that $\Mb$ is \textit{irreducible}, as we see from (\ref{ThreeParTransRate}) that the state space $\A$ constitutes a single communication class, i.e.,  for any  $\mb, \mb'\in\A$ we can find a finite path $\mb^{(1)}, \dots,\mb^{(l)}\in\A$ for which
$$
q(\mb^{(1)}|\mb)q(\mb^{(2)}|\mb^{(1)})\dots q(\mb'|\mb^{(l)})>0.
$$
Thus, we can argue as in Section 5 in \cite{Kendall1975} to draw some preliminary conclusions concerning the influence of the death rate on the dynamics of the model, based upon the properties of the sample size process $S$. In particular, in view of  the one-to-one correspondence
$$
S(t)=\sum_{i\ge1}iM_i(t)=0 \ \iff \Mb(t)=\eb_0=(0,0,\dots), \ t\ge0,
$$
we deduce that the initial condition $\eb_0$ and, by irreducibility, the entire state space $\A$ itself, are characterized by the same recurrence and transience properties of the state $0$ for $S$.  Accordingly, we conclude that  $\Mb$ is transient for $0<\mu<1$, null recurrent for $\mu=1$ and ergodic for $\mu>1$.

%

%

We shall henceforth focus on the latter case, wherein the underlying population almost surely undergoes a sequence of extinctions, separated by intervals of times of finite expected length during which its size and family configuration fluctuate according to the description in Section \ref{Sec:DefMod}. Such dynamics were shown in Lemma \ref{Lem:SReversibility} to generate a reversible population size process. The following theorem, which states our main result, shows that the same property holds for  $\Mb$.

\begin{theorem}\label{Theo:MReversibility}

Let $\Mb$ be as in Definition \ref{Def:ThreeParProc}, with $0<\alpha<1$ and $\mu>1$. Then, $\Mb$ is reversible with respect to the distribution $\pi$ with weights given by (\ref{ThreeParRevDistr}).

\proof

That $\pi$ is a well-defined probability distribution on $\A$, as well as the explicit expression of the constant $C$ in (\ref{CNormConst}), will follow from the mixture representation below. We then proceed to show the detailed balance condition
\eq{\label{MBalCond}
\pi(\mb')q(\mb|\mb')=\pi(\mb)q(\mb'|\mb),\ \mb\neq \mb'.
}
which, since $\Mb$ is regular and irreducible, is equivalent to the reversibility.

In view of the expression for the transition rates (\ref{ThreeParTransRate}) of $\Mb$, we only need to show (\ref{MBalCond}) for the cases in which, for given $\mb\in \A$, we have $\mb'=\mb+\eb_1$, $\mb'=\mb-\eb_i+\eb_{i+1}$ for $i\ge1$, or finally $\mb'=\mb+\eb_i-\eb_{i+1}$ for some $i\ge0$. The rate associated to the other transitions are indeed null.

For $\mb\in \A$ and $\mb'=\mb+\eb_1$, denoting $k=k(\mb),\ s=s(\mb),$ so that $k(\mb+\eb_1)=k+1$ and $s(\mb+\eb_1)=s+1$, the left hand side of (\ref{MBalCond}) becomes
\begin{equation}\nonumber
\begin{aligned}
\pi(\mb+\eb_1)&\,q(\mb|\mb+\eb_1)\\
=&\,C\bigg(\frac{\theta}{\alpha}\bigg)_{(k+1)}\textrm{Po}(m_1+1;\alpha_1\mu^{-1})\prod_{i\ge2}\mathrm{Po}(m_i;\alpha_i\mu^{-i})\mu(m_1+1)\\
=&\,
C\bigg(\frac{\theta}{\alpha}\bigg)_{(k)}\bigg(\frac{\theta}{\alpha}+k\bigg)\textrm{Po}(m_1;\alpha_1\mu^{-1})\frac{\alpha_1\mu^{-1}}{m_1+1}\prod_{i\ge2}\mathrm{Po}(m_i;\alpha_i\mu^{-i})\mu(m_1+1)\\
=&\,
C\bigg(\frac{\theta}{\alpha}\bigg)_{(k)}\prod_{i\ge1}\mathrm{Po}(m_i;\alpha_i\mu^{-i})(\theta+\alpha k)\\
=&\,
\pi(\mb)q(\mb+\eb_1|\mb).
\end{aligned}
\end{equation} 
Here we have used the identity
$$
\textrm{Po}(m_1+1;\alpha_1\mu^{-1})=e^{-\alpha_1\mu^{-1}}\frac{(\alpha_1\mu^{-1})^{m_1+1}}{(m_1+1)!}=\textrm{Po}(m_1;\alpha_1\mu^{-1})\frac{\alpha_1\mu^{-1}}{m_1+1}
$$
and the fact that
$$
\alpha_i=\frac{\alpha(1-\alpha)_{(i-1)}}{i!},\ i\ge1
$$
 with the convention $(1-\alpha)_{(0)}=1$, yielding $\alpha_1=\alpha$.

Instead, for $\mb'=\mb-\eb_i+\eb_{i+1}$, denoting as before $k=k(\mb)$ and $ s=s(\mb),$ we have $k(\mb-\eb_i+\eb_{i+1})=k$ and $s(\mb-\eb_i+\eb_{i+1})=s+1$. Thus, proceeding similarly
\eqstar{
\pi&(\mb-\eb_i+\eb_{i+1})q(\mb|\mb-\eb_i+\eb_{i+1})\\
&=
C\bigg(\frac{\theta}{\alpha}\bigg)_{(k)}
\prod_{i\ge1}\mathrm{Po}(m_i;\alpha_i\mu^{-i})
\dfrac{\alpha_{i+1}\mu^{-(i+1)}}{m_{i+1}+1}
\dfrac{m_i}{\alpha_i\mu^{-i}}\mu(i+1)(m_{i+1}+1)\\
&=\pi(\mb)\dfrac{\alpha_{i+1}}{\alpha_i}(i+1)m_i
=\pi(\mb)(i-\alpha)m_i\\
&=\pi(\mb)q(\mb-\eb_i+\eb_{i+1}|\mb),
}
where we have used
$$
\dfrac{\alpha_{i+1}}{\alpha_i}=\frac{i-\alpha}{i+1}.
$$

Finally, consider the case $\mb'=\mb+\eb_i-\eb_{i+1}$. Using the same notation as above, we have $k(\mb+\eb_i-\eb_{i+1})=k$ and $s(\mb+\eb_i-\eb_{i+1})=s-1$. A similar computation yields 
\eqstar{
\pi&(\mb+\eb_i-\eb_{i+1})q(\mb|\mb+\eb_i-\eb_{i+1})\\
&=
C\bigg(\frac{\theta}{\alpha}\bigg)_{(k)}
\prod_{i\ge1}\mathrm{Po}(m_i;\alpha_i\mu^{-i})
\dfrac{m_{i+1}}{\alpha_{i+1}\mu^{-(i+1)}}
\dfrac{\alpha_i\mu^{-i}}{m_i+1}(i-\alpha)(m_i+1)\\
&=\pi(\mb)\dfrac{\alpha_i}{\alpha_{i+1}}(i-\alpha)\mu m_{i+1}\\
&=\pi(\mb)\mu m_{i+1}(i+1)\\
&=\pi(\mb)q(\mb+\eb_i-\eb_{i+1}|\mb).
}
concluding the proof.
\endproof
\end{theorem}


The reversible distribution $\pi$ in of Theorem \ref{Theo:MReversibility} can be represented as a mixture of \PSFs  in (\ref{PSF}) with respect to the reversible $\textrm{N-Bin}(\theta,\mu^{-1})$ distribution of the population size process $S$. Specifically, we have
\eq{\label{MixRepr}
\pi(\mb)=\sum_{n\ge0}\mathrm{PSF}^{\alpha,\theta}_n(\mb)\textrm{N-Bin}(n;\theta,\mu^{-1}), \ \mb\in\A.
}
Indeed, for $\alpha>0$,  the PSF (\ref{PSF}) can be rewritten as
$$
\mathrm{PSF}^{\alpha,\theta}_n(\mb)=\dfrac{n!(\theta/\alpha)_{(k(\mb))}}{\theta_{(n)}}\prod_{i\ge1}\dfrac{\alpha_i^{m_i}}{m_i!}\1_{\{s(\mb)=n\}}, \ \mb\in\A
$$
where $\alpha_i$ is as in (\ref{ThreeParRevDistr}). Then, for $\mb\in\A$, it follows that
\begin{equation}\nonumber
\begin{aligned}
\sum_{n\ge0}&\,\mathrm{PSF}^{\alpha,\theta}_n(\mb)\textrm{N-Bin}(n;\theta,\mu^{-1})\\
=&\,
\sum_{n\ge0}\dfrac{n!(\theta/\alpha)_{(k(\mb))}}{\theta_{(n)}}\prod_{i\ge1}\dfrac{\alpha_i^{m_i}}{m_i!}\1_{\{s(\mb)=n\}}\frac{\theta_{(n)}}{n!}\mu^{-n}(1-1/\mu)^{\theta}\\
=&\,
(\theta/\alpha)_{(k(\mb))}\prod_{i\ge1}\dfrac{\alpha_i^{m_i}}{m_i!}\mu^{-s(\mb)}(1-1/\mu)^{\theta}\\
=&\,
(1-1/\mu)^{\theta}(\theta/\alpha)_{(k(\mb))}\prod_{i\ge1}\dfrac{\alpha_i^{m_i}}{m_i!}\mu^{-\sum_{i\ge1}im_i}\\
=&\,
(1-1/\mu)^{\theta}(\theta/\alpha)_{(k(\mb))}\prod_{i\ge1}\dfrac{(\alpha_i\mu^{-i})^{m_i}}{m_i!}e^{-\alpha_i\mu^{-i}}e^{\alpha_i\mu^{-i}}\\
=&\,(1-1/\mu)^\theta e^{\sum_{i\ge1}\alpha_i\mu^{-i}}(\theta/\alpha)_{(k(\mb))}\prod_{i\ge1}\mathrm{Po}(m_i;\alpha_i\mu^{-i}).
\end{aligned}
\end{equation} 
The series appearing as the argument of the exponential can be simplified by noticing that
\eqstar{
1-(1&-1/\mu)^{\alpha}=1-\sum_{i\ge0}\binom{\alpha}{i}(-\mu^{-i})\\
&=1-1-\sum_{i\ge1}(-1)^n\dfrac{\alpha(\alpha-1)\dots(\alpha-i+1)}{i!}\mu^{-i}\\
&=-\sum_{i\ge1}-\dfrac{-\alpha(1-\alpha)\dots(1-\alpha+i-2)}{i!}\mu^{-i}\\
&=\sum_{i\ge1}\dfrac{\alpha(1-\alpha)_{(i-1)}}{i!}\mu^{-i}
=\sum_{i\ge1}\alpha_i\mu^{-i}.
}
Replaced into the last line of the previous display, this recovers the expression for the reversible distribution $\pi$ with the constant $C$ as in (\ref{CNormConst}), implying (\ref{MixRepr}).


Theorem \ref{Theo:MReversibility} yields interesting insights in the long-run behavior of the three-parameter process: in particular, in view of the ergodicity, $\Mb(t)$ will converge in distribution to a stationary random allelic partition whose distribution admits the mixture representation (\ref{MixRepr}), and is characterized by a random, but almost surely finite number of items and groups. Thus, the reversible regime of $\Mb$ reveals clear analogies with the results in \cite{Kendall1975} for the BDI process, essentially due to the equivalence of the induced population size processes. However, as opposed to the case $\alpha=0$, $\Mb$  in the limit to maintains the mutual dependence among the multiplicities, when $\alpha>0$.

Concerning the range $0<\mu\le1$, we deduce at once from the absence of positive recurrence that the process is necessarily non-stationary, and in fact, because of the slower death rate, the associated population size process $S$ is seen in (\ref{PopSizeMargDistr}) to diverge as $t\to\infty$. On the other hand, for $\theta>0$, we can deduce from the description of the population model in Section \ref{Sec:DefMod}, that the number of families $K(t)$ alive at time $t$ stochastically dominates the number of families in the corresponding model with $\alpha=0$, as it can be rewritten as
$$
K(t)=K_1(t)+K_2(t)
$$
where $K_1(t)$ counts the family founded by immigrants, and $K_2(t)$ is the number of those started by the newborns leaving their family. Now $K_1(t)$ is evidently equal in distribution to the number of families comprised in a population governed by the BDI process, which, in turn, exhibits the logarithmic growth (\ref{KnConv1}). Accordingly, we conclude that $\Mb$ will evolve towards infinite structures characterized by an infinite number of individuals and families. Whether a precise description of the limiting behavior similar to (\ref{TwoParConv}) can be achieved within this formulation, remains  for the moment an open question.  


\section*{Acknowledgements}

The authors are grateful to the Associate Editor and two referees for helpful comments. The first author is supported by the ERC grant No. 647812 (UQMSI) and partially by the EPSRC grant EP/L016516/1 for the Cambridge Centre for Analysis.
The second and third authors are partially supported by the Italian Ministry of Education, University and Research (MIUR), through PRIN 2015SNS29B and through ``Dipartimenti di Eccellenza'' grant 2018-2022.

\bibliographystyle{agsm}
\bibliography{PopGenRef}

\end{document}